\documentclass[10pt]{article} 

\usepackage{amsmath}
\usepackage{amssymb} 
\usepackage{latexsym} 
\usepackage{epic,eepic,epsfig,psfrag} 
\usepackage{amscd}  
\usepackage{index}  

\newcommand{\QED}{\hspace*{\fill}$\Box$\medskip} 
\def\one{\hbox{1\hskip-2.7pt l}}
\def\smone{{\scriptstyle\rm 1\hskip-2.05pt l}}
\def\ssmone{{\scriptscriptstyle\rm 1\hskip-2.05pt l}}

\def\rd{{\rm d}}

\def\p{\phi}

\def\d{\delta} 

\def\ep{\varepsilon} 
\def\e{\eta}

\def\r{\rho}

\def\n{\nu} 
\def\m{\mu}

\def\G{\Gamma} 

\def\S{\Sigma}

\def\cC{{\mathcal C}}

\def\R{{\mathbb R}} 
  
\def\N{{\mathbb N}} 
 
\def\H{{\mathbb H}}

\def\half{{\textstyle{\frac 12}}} 

\def\laplace{\Delta} 
\def\Pr{{\bf Proof:}\;} 
\def\st{\: \big| \:}

\def\dt{{\rm d}t}

\def\dr{{\rm d}r}
\def\pd{\partial}

\def\bc{{\bar c}}

\def\bx{{\bar x}}
\def\by{{\bar y}}
\def\brh{{\bar\r}}

\def\dvol{\,\rd{\rm vol}}

\def\tint{\textstyle\int}

\newtheorem{dfn}{Definition}[section] 
\newtheorem{lem}[dfn]{Lemma} 
\newtheorem{prp}[dfn]{Proposition} 
\newtheorem{thm}[dfn]{Theorem} 
\newtheorem{rmk}[dfn]{Remark}

\begin{document}

\bibliographystyle{plain}

\author{Katrin Wehrheim 
}

\title{Energy quantization and mean value inequalities for nonlinear boundary value
       problems}

\maketitle

\begin{abstract}
We give a unified statement and proof of a class of wellknown mean value inequalities 
for nonnegative functions with a nonlinear bound on the Laplacian. We generalize these
to domains with boundary, requiring a (possibly nonlinear) bound on the normal derivative
at the boundary.
These inequalities give rise to an energy quantization principle for sequences of solutions of
boundary value problems that have bounded energy and whose energy densities satisfy nonlinear
bounds on the Laplacian and normal derivative: One obtains local uniform bounds on the complement
of finitely many points, where some minimum quantum of energy concentrates.
\end{abstract}

\section{Energy quantization}
\label{eq}

One purpose of this note is to explain an 'energy quantization' principle that 
-- in different forms -- has successfully been applied to a variety of 
partial differential equations, such as minimal submanifolds, harmonic maps, 
pseudoholomorphic curves, and Yang-Mills connections.
The common feature of these PDE's is an energy functional. (The solutions often but
not necessarily are critical points thereof.)
The phenomenon that we call 'energy quantization' is a consequence of a mean value
inequality for the energy density.
The second purpose of this note is to give a unified statement and proof of the underlying 
mean value inequality for the Laplace operator.
In section~\ref{mean} we moreover generalize this inequality to domains with boundary and
inhomogeneous Neumann boundary conditions.

In this section we generally consider a PDE for maps $u:D\to T$ from a Riemannian manifold $D$
(with possibly nonempty boundary $\pd D$)
to a target space $T$, e.g.\ another manifold, a Banach space, or a fibre bundle over $D$.
The energy is given for all sufficiently regular maps $u$ in the form
$$
E(u) = \int_D e(u) ,
$$
where the integrand $e(u):D\to[0,\infty)$ is a nonnegative energy density function. 
Its key property is that for a solution $u$ of the PDE, the positive definite Laplacian 
$\laplace e(u)$ can be bounded above in terms of $e(u)$ itself.

If this bound is linear in $e(u)$, then one immediately obtains a $\cC^0$-control on
$e(u)$ in terms of its mean values on geodesic balls. If moreover the energy $E(u_i)$ of a 
sequence of solutions $u_i$ is bounded, then this provides a uniform bound on the energy 
densities $e(u_i)$ on any compact subdomain of $D\setminus\pd D$.
In many cases this leads to a compactness result, i.e.\ to the convergence of a subsequence
of the solutions $u_i$.

For solutions of nonlinear PDE's however, the bound on $\laplace e(u)$ is usually nonlinear in $e(u)$.
In that case, for polynomial nonlinearities up to some order depending on the $\dim D$,
the mean value inequality only hold on geodesic balls with sufficiently small energy.
In theorem~\ref{thm mean} we give a precise and general statement of this wellknown fact.
In theorem~\ref{thm mean bdy}, we then generalize this mean value inequality to domains $D$ 
with boundary $\pd D$ and bounds on the outer normal derivative $\frac\pd{\pd\n} e|_{\pd D}$.
This allows to obtain uniform bounds on $e(u_i)$ up to the boundary for solutions $u_i$ of
a PDE with appropriate boundary conditions and with bounded energy.

In the case of nonlinear bounds on the Laplacian or the normal derivative one only obtains
locally uniform bounds on the complement of finitely many points:
By a converse of the mean value inequality, the energy densities $e(u_i)$ can only blow up at 
points where some nonzero quantum of energy concentrates.
In the following lemma we give a blueprint for such energy quantization results.

Here $D$ is a Riemannian manifold (possibly noncompact or with boundary),
$\laplace=\rd^*\rd$ denotes the positive definite Laplace operator, and
$\frac\pd{\pd\n}$ denotes the outer normal derivative at $\pd D$.
For the sake of simplicity we make a technical assumption on the metric near the boundary.

\medskip

\noindent
{\bf Assumption:} A neighbourhood of $\pd D\subset D$ 
is locally isometric to the Euclidean half space $\H^n$. 

\medskip

For a general metric near the boundary, the mean value inequality becomes harder and more technical to
prove, but theorem~\ref{thm mean bdy} should generalize in the same way as theorem~\ref{thm mean},
and hence this lemma should extend to general Riemannian manifolds with boundary.

\begin{lem}\label{lem eq}
There exists a constant $\hbar>0$ depending on $n=\dim D$ and some constants $a,b\geq 0$ 
such that the following holds:
Let $e_i\in\cC^2(D,[0,\infty))$ be a sequence of nonnegative functions such that for some constants
\hbox{$A_0,A_1,B_0,B_1\geq 0$}
\begin{equation*}
\left\{\begin{array}{ll}
\laplace e_i &\leq  A_0 + A_1 e_i + a e_i^{\frac {n+2}n} , \\
\frac\pd{\pd\n}\bigr|_{\pd D} e_i \hspace{-2mm} &\leq B_0 + B_1 e_i + b e_i^{\frac {n+1}n} .
\end{array}\right.
\end{equation*}
Moreover, suppose that there is a uniform bound\; $\displaystyle \int_D e_i \leq E <\infty$.

Then there exist finitely many points, $x_1,\ldots,x_N\in D$ (with $N\leq E/\hbar$) and a subsequence 
such that the $e_i$ are uniformly bounded on every compact subset of $D\setminus\{x_1,\ldots x_N\}$, 
and there is a concentration of energy $\hbar>0$ at each $x_j$ :
For every $\d>0$ there exists $I_{j,\d}\in\N$ such that
\begin{equation}\label{concentration}
 \int_{B_\d(x_j)} e_i \geq \hbar \qquad\forall i\geq I_{j,\d}.
\end{equation}
\end{lem}
\Pr
Suppose that for some $x\in D$ there is no neighbourhood on which the $e_i$ are uniformly bounded.
Then there exists a subsequence (again denoted $(e_i)$) and $D\ni z_i\to x$ such that
$e_i(z_i)=R_i^n$ with $R_i\to\infty$.
We can then apply the mean value inequality theorem~\ref{thm mean} on the balls $B_{\d_i}(z_i)$
of radius $\d_i=R_i^{-\frac 12}>0$. For sufficiently large $i\in\N$, these lie within appropriate 
coordinate charts of $D$. In case $z\in\pd D$ we use the Euclidean coordinate charts at the boundary 
to apply theorem~\ref{thm mean bdy}, but we also denote the balls in half space by $B_{\d_i}(z_i)$.
These mean value theorems provide uniform constants $C$ and $\hbar:=\m(a,b)>0$ such 
that for every $i\in\N$ either
\begin{equation}\label{hbar}
\int_{B_{\d_i}(z_i)} e_i > \hbar
\end{equation}
or $\int_{B_{\d_i}(z_i)} e_i \leq \hbar$ and hence
$$
R_i^n \;=\; e(z_i) 
\;\leq\; C A_0 \d_i^2 + C B_0 \d_i +  C \bigl(A_1^{\frac n2} + B_1^n + \d_i^{-n} \bigr)\hbar .
$$
In the latter case multiplication by $\d_i^n=R_i^{-\frac n2}$ implies
\begin{equation}\label{eq est}
R_i^{\frac n2} \;\leq\; C A_0 R_i^{-\frac{n+2}2} + C B_0 R_i^{-\frac{n+1}2} 
                   +  C\hbar \bigl(A_1^{\frac n2}R_i^{-\frac n2} + B_1^nR_i^{-\frac n2} + 1 \bigr).
\end{equation}
As $i\to\infty$, the left hand side diverges to $\infty$, whereas the right hand side converges
to $C\hbar$. Thus the alternative (\ref{hbar}) must hold for all sufficiently large $i\in\N$.
In particular, this implies the energy concentration (\ref{concentration}) at $x_j=x$.

Now we can go through the same argument for any other point $x\in D$ at which the present 
subsequence $(e_i)$ is not locally uniformly bounded.
That way we iteratively find points $x_j\in D$ such that the energy concentration 
(\ref{concentration}) holds for a further subsequence $(e_i)$. 
Suppose this iteration yields $N>E/\hbar$ distinct points $x_1,\ldots,x_N$ (and might not even
terminate after that). Then we would have a subsequence $(e_i)$ for which at least energy $\hbar>0$
concentrates near each $x_j$. Since the points are distinct, this contradicts the energy bound
$\int_D e_i \leq E$.
Hence this iteration must stop after at most $\lfloor E/\hbar \rfloor$ steps, 
when the present subsequence $(e_i)$ is locally uniformly bounded in the complement of
the finitely many points, where we found the energy concentration before.
\QED

The analogy in the compactness proofs for a variety of PDE's, including minimal submanifolds 
\cite{A,CS} and harmonic maps of surfaces \cite{SU}, has already been observed and listed 
by Wolfson \cite{Wo}. In all cases, the nonlinearities are exactly of the maximal order as
in lemma~\ref{lem eq}.
Here we discuss pseudoholomorphic curves and Yang-Mills connections in more detail.

For {\bf pseudoholomorphic curves} (with a 2-dimensional domain) 
the energy is the $L^2$-norm of the gradient, and the estimate 
$\laplace e \leq C(e+e^2)$ leads to Gromov's compactness result \cite{G,Wo}.
For {\bf pseudoholomorphic curves with Lagrangian boundary conditions}, these compactness 
results are also wellknown. They can be proven via a specific choice of a metric for 
which $\frac\pd{\pd\n} e = 0$ (see \cite{F} and \cite[4.3]{MS}).
Then the energy density can be extended across the boundary by reflection and the mean value 
inequality for $\R^n$ applies.
For the metric given by the almost complex structure, the Lagrangian boundary condition only 
implies $\frac\pd{\pd\n} e \leq C(e+e^{\frac 32})$, which fits nicely into our energy quantization 
principle.

For {\bf Yang-Mills connections} on 4-manifolds the energy is the $L^2$-norm of the curvature.
The bound $\laplace e\leq C(e+e^{\frac 32})$ was used by Uhlenbeck \cite{U} to prove 
a removable singularity result, which leads to Donaldson's compactification of the 
moduli space of anti-self-dual instantons \cite{D}.
For a proof of the energy quantization as in lemma~\ref{lem eq} 
see also \cite[Thm.2.1]{W bubb}.

For {\bf anti-self-dual instantons with Lagrangian boundary conditions} on $\H^2\times\S$, 
an argument along the lines of lemma~\ref{lem eq} is used in \cite[Thm.1.2]{W bubb}.
There the energy density is given by the slicewise $L^2$-norm of the curvature, 
$e(A)=\int_\S |F_A|^2$.
The Lagrangian boundary condition (which has global nature along the Riemann surface $\S$) 
provides $\frac\pd{\pd\n} e \leq C(e+e^{\frac 32})$,
but one only has a linear bound $\laplace e \leq g e$ with a function $g$ that cannot be 
bounded in terms of $e$ or a constant. 
The argument using (\ref{eq est}) however allows for a mild blowup of $A_1$, and the according
estimate $|g|\leq C R_i^2$ can be established.
This result does not follow from the standard rescaling methods for Yang-Mills connections.

\section{Mean value inequalities}
\label{mean}

In this section we state and prove the mean value inequalities that were desribed in section~\ref{eq}
and that the energy quantization principle is based on.

The subsequent theorem is wellknown and proofs in an exhausting collection of cases can be found
in the literature, e.g.\ \cite{S,U}. Our aim here is to give a unified statement and proof.
In theorem~\ref{thm mean bdy} below we moreover generalize this result to domains with boundary.

We denote by $B_{r_0}(x_0)\subset\R^n$ the open geodesic ball of radius $r_0$ centered at 
$x_0\in\R^n$ and with respect to the present metric.
The Laplace operator $\laplace=\rd^*\rd$ and integration will also be using the metric given 
in the context.
The Euclidean metric on $\R^n$ is denoted by its matrix $\one$.

\begin{thm}
\label{thm mean}
For every $n\in\N$ there exist constants $C$, $\m>0$, and $\d>0$ such that the 
following holds for all metrics $g$ on $\R^n$ such that $\|g-{\emph\one}\|_{W^{1,\infty}}\leq\d$.

Let $B_r(0)\subset\R^n$ be the geodesic ball of radius $0<r\leq 1$. 
Suppose that $e\in\cC^2(B_r(0),[0,\infty))$ satisfies for some $A_0,A_1,a\geq 0$
$$
\laplace e \leq A_0 + A_1 e + a \, e^{\frac {n+2}n}
\qquad\text{and}\qquad
\int_{B_r(0)} e \leq \m a^{-\frac n2} .
$$
Then
$$
e(0) \leq  C A_0 r^2 + C \bigl( A_1^{\frac n2} + r^{-n} \bigr) \int_{B_r(0)} e .
$$
\end{thm}

\begin{rmk}
By using local geodesic coordinates the above theorem also implies a mean value inequality on 
closed Riemannian manifolds with uniform constants $C,\m>0$, and for all geodesic balls of radius 
less than a uniform constant.
\end{rmk}

The special case $A_0=A_1=a=0$ of theorem~\ref{thm mean} and the starting point for the proof is 
Morrey's \cite{M} mean value inequality for subharmonic functions. A proof of the version below can
be found in e.g.\ \cite[Thm.2.1]{LS}. For the Euclidean metric $g=\one$ we give an elementary
proof in lemma~\ref{lem 1} below.

\begin{prp} \label{prp subharmonic}
For every $n\in\N$ there exist constants $C$ and $\d>0$ such that the following 
holds for all $0<r\leq 1$ and all metrics $g$ on $\R^n$ with 
\hbox{$\|g-\emph\one\|_{W^{1,\infty}}\leq\d$}:
If $e\in\cC^2(B_r(0),[0,\infty))$ satisfies $\laplace e \leq 0$ then
$$
e(0) \leq C r^{-n} \int_{B_r(0)} e .
$$
\end{prp}

\noindent
{\bf Proof of theorem \ref{thm mean}: } \\
This proof is based on the Heinz trick, which is to consider the maximum $\bar c$ 
of the function $f$ below. This allows to replace the bound on the 
Laplacian by a constant depending on $\bar c$. One then obtains the result from the mean
value inequality for subharmonic functions and a number of rearrangements in different cases.

Consider the function $f(\r)=(1-\r)^n \sup_{B_{\r r}(0)} e$ for $\r\in[0,1]$. 
It attains its maximum at some $\bar\r<1$.
Let $\bar c=\sup_{B_{\brh r}(0)}e = e(\bx)$ and $\ep=\half(1-\brh)<\half$, then
$$
e(0) \;=\; f(0) \;\leq\; f(\brh) \;=\; 2^n\ep^n \bc .
$$
Moreover, we have for all $x\in B_{\ep r}(\bx)\subset B_r(y)$
$$
e(x) 
\;\leq\; \sup_{B_{(\brh+\ep)r}(0)} e 
\;=\; \bigl(1-\brh-\ep)^{-n} f(\brh+\ep)
\;\leq\; 2^n (1-\brh)^{-n} f(\brh)
\;=\; 2^n \bc ,
$$
and hence $\laplace e \leq A_0 + 2^n A_1 \bc + 2^{n+2} a \bc^{\frac{n+2}n}$.
Now define the function 
$$
v(x):= e(x) + \tfrac 1n \bigl( A_0 + 2^n \bc \bigl( A_1 + 4 a \bc^{\frac 2n}\bigr)\bigr) |x-\bx|^2 
$$ 
with the Euclidean norm $|x-\bx|$. It is nonnegative and subharmonic 
on $B_{\ep r}(\bx)$ if the metric is sufficiently $\cC^1$-close to $\one$. 
This is since $\laplace_{\smone}|x-\bx|^2=2n$ for the Euclidean metric
and $|x-\bx|\leq \ep r \leq 1$ is bounded, so $\laplace|x-\bx|\geq n$
whenever $\|g-\one\|_{W^{1,\infty}}\leq\d$ is sufficiently small.
The control of the metric also ensures that
the integral $\int_{B_{\r r}(\bx)} |x-\bx|^2$ is bounded by the following 
integral over the Euclidean ball $B^\smone_{2\r r}(\bx)$
$$
2 \int_{B^\ssmone_{2\r r}(\bx)} |x-\bx|^2 
\;=\; 2  \int_0^{2\r r} t^{n+1} {\rm Vol}\, S^{n-1} \; \dt
\;=\; C_n (\r r)^{n+2} .
$$
So we obtain from proposition~\ref{prp subharmonic} with a uniform constant $C$ for all $0<\r\leq\ep$
$$
\bc \;=\; v(\bx) \;\leq\; C(\r r)^{-n} \int_{B_{\r r}(\bx)} v
 \;\leq\;  C \bigl(A_0 + \bc \bigl(A_1 + 4 a \bc^{\frac 2n}\bigr)\bigr) (\r r)^2 
          \;+\; C (\r r)^{-n} \int_{B_{\r r}(\bx)} e  .         
$$
If $C A_0 (\ep r )^2  \leq \half \bc$, then we can drop $A_0$ from this inequality, just changing the
constant $C$. Otherwise $e(0)\leq \bc \leq C A_0 r^2$ already proves the assertion.

Next, if $C(A_1+4a\bc^{\frac 2n}) (\r r)^2  \leq \half$, then this implies 
$\bc \leq 2C(\r r)^{-n}\int_{B_r(0)} e$.
So if $C(A_1+4a\bc^{\frac 2n}) (\ep r)^2  \leq \half$ then $\r=\ep$ proves the assertion,
$$
e(0) \;\leq\; 2^n \ep^n \bc \;\leq\; 2^{n+1} C r^{-n} \tint_{B_r(0)} e .
$$
Otherwise we can choose $0<\r<\ep$ such that 
$\r r = \bigl( 2 C(A_1+4a\bc^{\frac 2n})\bigr)^{-\frac 12}$ 
to obtain with a new uniform constant $C$
$$
e(0) \;\leq\; \bc 
\;\leq\; C \bigl( A_1+4a\bc^{\frac 2n} \bigr)^{\frac n2}  \tint_{B_{\r r}(\bx)} e  .
$$
Again we have to distinguish two cases: Firstly, if $4a\bc^{\frac 2n}\leq A_1$ then this 
yields 
$$
e(0) \;\leq\; C (2A_1)^{\frac n2} \tint_{B_{\r r}(\bx)} e .
$$
Secondly, if $A_1<4a\bc^{\frac 2n}$ then 
$\bc < \bc\,C (8a)^{\frac n2} \int_{B_{\r r}(\bx)} e $ and thus with some $\m>0$
$$
\tint_{B_r(0)} e \;>\; \m a^{-\frac n2}.
$$
So altogether we either have the above or with some uniform constant $C$
$$
e(0) \;\leq\; C A_0 r^2 + C \bigl( A_1^{\frac n2} + r^{-n} \bigr) \tint_{B_r(0)} e .
$$
\QED

Next we generalize the mean value inequality theorem~\ref{thm mean} to the half space
$$
\H^n=\{(x_0,\bx) \st x_0 \in [0,\infty) , \bx \in \R^{n-1} \} .
$$
In order to generalize the mean value inequality to manifolds with boundary we would have to
consider general metrics on $\H^n$.
This would however disturb the elementary geometric proof, so we restrict this exposition to 
the case of the Euclidean metric.
In that case the outer normal derivative $\frac\pd{\pd\n}|_{\pd\H^n}$ is 
just $-\frac\pd{\pd x_0}|_{x_0=0}$.
We moreover denote the intersection of a Euclidean ball with the half space by
$$ 
D_{r_0}(x_0):= B_{r_0}(x_0)\cap\H^n .
$$

\begin{thm}
\label{thm mean bdy}
For every $n\geq 2$ there exists a constant $C$ and for all $a,b\geq 0$ there exists 
$\m(a,b)>0$ such that the following holds: 

Let $D_r(y)\subset\H^n$ be the Euclidean $n$-ball of radius $r>0$ and center 
$y\in\H^n$ intersected with the half space. 
Suppose that $e\in\cC^2(D_r(y),[0,\infty))$ satisfies for some $A_0,A_1,B_0,B_1\geq 0$
\begin{equation*}
\left\{\begin{array}{ll}
\laplace e &\leq  A_0 + A_1 e + a e^{\frac {n+2}n} , \\
\frac\pd{\pd\n}\bigr|_{\pd\H^n} e \hspace{-2mm} &\leq B_0 + B_1 e + b e^{\frac {n+1}n} ,
\end{array}\right.
\qquad\text{and}\qquad
\int_{D_r(y)} e \leq \m(a,b) .
\end{equation*}
Then
$$
e(y) \leq  C A_0 r^2 + C B_0 r +  C \bigl(A_1^{\frac n2} + B_1^n + r^{-n} \bigr) \int_{D_r(y)} e .
$$
\end{thm}

We will prove this theorem in three steps.
The first step is the generalization of proposition~\ref{prp subharmonic} to domains with 
boundary and subharmonic functions in the sense of the weak Neumann equation:
A distribution $e$ on a manifold $M$ is called subharmonic if for all $\psi\in\cC^\infty(M,[0,\infty))$
with $\frac{\pd\psi}{\pd\n}|_{\pd M}=0$
$$
0 \;\geq\; \int_M e \, \laplace\psi 
\;=\; \int_M \psi\,\laplace e + \int_{\pd M} \psi\,\tfrac{\pd e}{\pd\n} .
$$
For $e\in\cC^2(M)$ the equality above holds and implies that $\laplace e \leq 0$ and
$\frac{\pd e}{\pd\n}|_{\pd M}\leq 0$.

\begin{lem} \label{lem 1}
There exists a constant $C$ such that the following holds for all $R>0$ and $y\in\H^n$:
Suppose that $e\in\cC^2(D_r(y),[0,\infty))$ satisfies
\begin{equation*}
\left\{\begin{array}{ll}
\laplace e &\leq 0 , \\
\frac\pd{\pd\n}\bigr|_{\pd\H^n} e \hspace{-2mm} &\leq 0 .
\end{array}\right.
\end{equation*}
Then
$$
e(y) \leq C r^{-n} \int_{D_r(y)} e .
$$
\end{lem}
\Pr
We write $\H^n=\{(x_0,\bx) \st x_0 \in [0,\infty) , \bx \in \R^{n-1} \}$
and use spherical coordinates
$(x_0,\bx) = ( y_0 + r \cos\p \,,\, \by + r\sin\p \cdot z ) =: (r,\p,z)$ 
with $r\in[0,\infty)$, $\p\in[0,\pi]$, and $z\in S^{n-2}\subset\R^n$.
(For $n=2$ this notation means $S^0=\{-1,1\}$, and integration $\int_{S^0} \ldots \dvol_{S^0}$
will denote summation of the values at these two points.)
Now the boundary of $D_r(y)$ has two parts,
\begin{align*}
Z_r &:=\; \pd D_r(y) \cap \pd\H^n 
\;=\; \bigl\{ (0 \,,\, \bx ) \st |\bx-\by|^2 \leq r^2 - {y_0}^2 \bigr\} , \\
\G_r &:=\;\pd B_r(y) \cap \H^n
\;=\; \bigl\{ (r,\p,z) \st \p\in[0,\p_0(r)], z\in S^{n-2} \bigr\} .
\end{align*}
Here we use $\p_0(r):=\arccos(-y_0/r)$. 
For $y_0>r$ we set $\p_0(r):=\pi$, so $\G_r$ is the entire sphere and the set $Z_r$ 
is empty. With this we calculate for all $r>0$
\begin{align}
&\frac\rd\dr \left( r^{-n+1} \int_{\G_r} e \right) \nonumber\\
&= \frac\rd\dr \left( r^{-n+1} \int_0^{\p_0(r)}\int_{S^{n-2}}
                      e(r,\p,z) \,(r\sin\p)^{n-2} \dvol_{S^{n-2}} \, r \,\rd\p \right) 
\nonumber \\
&= \int_0^{\p_0(r)}\int_{S^{n-2}} 
      \pd_r e(r,\p,z) \,(\sin\p)^{n-2} \dvol_{S^{n-2}} \, \rd\p   \label{dr int}\\
&\quad + \frac{\pd\p_0}{\pd r} \int_{S^{n-2}} 
       e(r,\p_0(r),z) \,(\sin\p_0(r))^{n-2} \dvol_{S^{n-2}} .  \nonumber
\end{align}
Note that $\p_0(r)$ is constant for $y_0=0$ as well as for $r\leq y_0$. So firstly
in case $y_0>0$ we have for all $0<r\leq y_0$
\begin{equation} \label{int small r}
\frac\rd\dr \left( r^{-n+1} \int_{\G_r} e \right) 
\;=\; r^{-n+1} \int_{\pd D_r(y)} \tfrac\pd{\pd\n} e  
\;=\; - r^{-n+1} \int_{D_r(y)} \laplace e  
\;\geq\; 0 .
\end{equation}
In that case we moreover have 
\begin{equation}\label{lim int}
\lim_{r\to 0} \left( r^{-n+1} \int_{\G_r} e \right) \;=\; {\rm Vol}\,S^{n-1} \, e(y) ,
\end{equation}
so integrating $\int_0^{\frac R2} r^{n-1} \ldots \,\dr$ proves the lemma for all 
$R\leq 2 y_0$,
$$
\tfrac 1n 2^{-n} R^n \,{\rm Vol}\,S^{n-1} \, e(y)
\;\leq\; \int_0^{\frac R2} \int_{\G_r} e \, \dr 
\;\leq\; \int_{D_R(y)} e .
$$
Next, in case $y_0=0$ we have for all $r>0$
\begin{align*}
\frac\rd\dr \left( r^{-n+1} \int_{\G_r} e \right) 
&=\; r^{-n+1} \int_{\G_r} \tfrac\pd{\pd\n} e  \\
&=\; - r^{-n+1} \int_{D_r(y)} \laplace e 
\; - \; r^{-n+1} \int_{Z_r} \tfrac\pd{\pd\n} e  
\;\geq\; 0 .
\end{align*}
Since $\lim_{r\to 0} \left( r^{-n+1} \int_{\G_r} e \right) 
= \half {\rm Vol}\,S^{n-1} \, e(y)$, integration over $0<r\leq R$ then proves
the lemma for $y_0=0$ and all $R>0$,
$$
\tfrac 1{2n} R^n \,{\rm Vol}\,S^{n-1} \, e(y)
\;\leq\; \int_0^R \int_{\G_r} e \, \dr 
\;=\; \int_{D_R(y)} e .
$$
Finally, in case $R>2 y_0>0$ we obtain from (\ref{dr int}) for all $r > y_0$
\begin{align*}
&\frac\rd\dr \left( r^{-n+1} \int_{\G_r} e \right) 
&\geq 
\frac{- y_0}{r\sqrt{r^2-y_0^2}} \int_{S^{n-2}} 
  e(r,\p_0(r),z)  \,(\sin\p_0(r))^{n-2}\dvol_{S^{n-2}} . 
\end{align*}
Now we can use (\ref{lim int}), (\ref{int small r}), and integrate the above to obtain
for all $y_0 < r \leq \half R$
\begin{align*}
{\rm Vol}\,S^{n-1} \, e(y) 
&\leq r^{-n+1} \int_{\G_r} e \\
&\quad + \int_{y_0}^r y_0 \,\r^{1-n} (\r^2-y_0^2)^{\frac{n-3}2} \int_{S^{n-2}} 
                            e(\r,\p_0(\r),z) \, \dvol_{S^{n-2}}\,\rd\r . 
\end{align*}
Since $(\r,\p_0(\r),z)\in\pd\H^n$, we already know that
$$
e(\r,\p_0(\r),z) 
\;\leq\; \frac{2n}{{\rm Vol}\,S^{n-1}(\frac R2)^n } \int_{D_{\frac R2}(\r,\p_0(\r),z)} e
\;\leq\; \frac{2^{n+1}n}{{\rm Vol}\,S^{n-1} R^n} \int_{D_R(y)} e .
$$
With this we find that for all $0 < y_0 < r \leq \half R$
\begin{align}
{\rm Vol}\,S^{n-1} \, e(y)
&\leq r^{-n+1} \int_{\G_r} e \nonumber \\
&\quad +
\frac{2^{n+1}n\, {\rm Vol}\,S^{n-2}}{R^n\, {\rm Vol}\,S^{n-1}} 
\int_1^{r y_0^{-1}} t^{-2} \bigl( 1 - t^{-2} \bigr)^{\frac{n-3}2}\,\dt 
\; \int_{D_R(y)} e \nonumber\\
&\leq r^{-n+1} \int_{\G_r} e \;+\; C_n R^{-n}\int_{D_R(y)} e .  \label{int large r}
\end{align}
Here we have introduced a constant $C_n$ that only depends on $n\geq 2$, in particular
on the value of the integral in $t$: For $n=2$ we calculate it explicitly,
$$
\int_1^{r y_0^{-1}} t^{-2} \bigl( 1 - t^{-2} \bigr)^{-\frac 12}\,\dt
\;=\; \Bigl[ \arccos(t^{-1}) \Bigr]_1^{r y_0^{-1}}
\;=\; \arccos\bigl(\tfrac r {y_0}\bigr)
\;<\; \tfrac \pi 2 .
$$
For $n\geq 3$ we have
$$
\int_1^{r y_0^{-1}} t^{-2} \bigl( 1 - t^{-2} \bigr)^{\frac{n-3}2}\,\dt
\;\leq \; \int_1^{r y_0^{-1}} t^{-2} \dt
\;=\; 1 - \tfrac r {y_0}
\;<\; 1 .
$$
Now from (\ref{int small r}) we know that (\ref{int large r}) also holds for
$0<r\leq y_0$ (with $C_n=0$), so integrating $\int_0^{\frac R2} r^{n-1} \ldots \,\dr$ 
proves the lemma in this last case,
$$
\tfrac 1n \bigl(\tfrac R2\bigr)^n {\rm Vol}\,S^{n-1} \, e(y)
\;\leq\; \int_0^{\frac R2} \int_{\G_r} e \,\dr 
\;+\; \tfrac 1n \bigl(\tfrac R2\bigr)^n C_n R^{-n}\int_{D_R(y)} e \\
\;\leq\; C \int_{D_R(y)} e .
$$

\vspace{-2.5mm}

\hspace*{\fill}$\Box$
\medskip

\noindent
{\bf Proof of theorem \ref{thm mean bdy}: } \\
With lemma~\ref{lem 1} in hand, the second step of the proof is to assume constant positive
bounds, $\laplace e \leq A$ and $\frac\pd{\pd\n}\bigr|_{\pd\H^n} e \leq B$
and find a constant $C$ such that for all $r>0$ and $y\in\H^n$
\begin{equation} \label{A B bounds}
e(y) \;\leq\; C r^{-n} \int_{D_r(y)} e \;+\; C A r^2 + C B r.
\end{equation}
That is, we first prove the theorem with $A_1=B_1=a=b=0$. To do this consider the function
$$
v(x):= e(x) + \tfrac 1{2n} A |x-y|^2 + ( B + \tfrac 1n A \, y_0 ) x_0 .
$$
It is positive and satisfies $\laplace v \leq 0$ and
$\tfrac\pd{\pd\n}\bigr|_{\pd\H^n} v \leq 0$, so lemma~\ref{lem 1} implies that
\begin{equation}  \label{step 1 calc}
e(y) 
\;=\; v(y) - ( B + \tfrac 1n A \, y_0 ) y_0 
\;\leq\; v(y) 
\;\leq\; C r^{-n} \int_{D_r(y)} v .
\end{equation}
In case $r\leq y_0$ we just use $v(x)= e(x) + \tfrac 1{2n} A |x-y|^2$, then the
same holds, and moreover 
$$
\int_{D_r(y)} v 
\;=\; \int_{D_r(y)} e \;+\; \tfrac 1{2n} A \int_0^r t^{n+1}{\rm Vol}\,S^{n-1}\dt
\;=\; \int_{D_r(y)} e \;+\; \tfrac {{\rm Vol}\,S^{n-1}}{2n(n+2)} A r^{n+2} .
$$
In case $r>y_0$ we have (using $x_0\leq 2r$ on $B_r(y)$)
\begin{align*}
\int_{D_r(y)} v 
&\leq \int_{D_r(y)} e \;+\; \tfrac 1{2n} A \int_0^r t^{n+1}{\rm Vol}\,S^{n-1}\dt 
\;+\; ( B + \tfrac 1n A \, y_0 ) \int_{D_r(y)} x_0  \\
&\leq \int_{D_r(y)} e \;+\; \tfrac {{\rm Vol}\,S^{n-1}}{2n(n+2)} A r^{n+2} 
\;+\; ( B + \tfrac 1n A \, r ) \tfrac 2n {\rm Vol}\,S^{n-1} r^{n+1} .
\end{align*}
In any case, putting this into (\ref{step 1 calc}) proves (\ref{A B bounds}).

\pagebreak

Finally, to prove the theorem we consider
-- analogous to the proof of theorem~\ref{thm mean} -- the function 
$f(\r)=(1-\r)^n \sup_{D_{\r r}(y)} e$ defined for $\r\in[0,1]$. 
It attains its maximum at some $\bar\r<1$.
We denote $\bar c=\sup_{D_{\brh r}(y)}e = e(\bx)$ and $\ep=\half(1-\brh)$, then
$e(y) \leq 2^n\ep^n \bc$ and $e(x) \leq 2^n \bc$ for all $x\in D_{\ep r}(\bx)$.
Thus on $D_{\ep r}(\bx)\subset D_r(y)$ we have 
$\laplace e \,\leq\, A_0 + 2^n\bc ( A_1 + 4a\bc^{\frac 2n}) $
and $\tfrac\pd{\pd\n}\bigr|_{\pd\H^n} e \leq B_1 + 2^n\bc (B_1 + 2b\bc^{\frac 1n})$.
Putting this into (\ref{A B bounds}) yields for all $0<\r\leq\ep$
\begin{align}
\bc \;=\; e(\bx) 
&\;\leq\; C (\r r)^{-n} \tint_{D_{\r r}(\bx)} e 
\;+ C \bigl(A_0 + 2^n\bc ( A_1 + 4a\bc^{\frac 2n})\bigr) (\r r)^2 \label{basic}\\
&\qquad\qquad\qquad\qquad\quad\;\;\;
+ C \bigl(B_0 + 2^n\bc (B_1 + 2b\bc^{\frac 1n})\bigr) \r r .  \nonumber
\end{align}
To deduce the claimed mean value inequality from this, we have to go through a number
of different cases.
Firstly, if $CA_0(\ep r)^2 + CB_0\ep r \geq \half \bc$, then since $\ep\leq\half$
$$
e(y) \;\leq\; \bc \;\leq \; C A_0 r^2 + C B_0 r ,
$$
which proves the theorem. Otherwise (\ref{basic}) continues to hold with $A_0$ and $B_0$
dropped (and another constant).
Next, let $0<\ep'<\ep$ be the solution of $ A_1(\ep' r)^2 + B_1\ep' r = 2^{-n-1} C^{-1}$ 
or in case $A_1(\ep r)^2 + B_1\ep r \leq 2^{-n-1} C^{-1}$ let $\ep=\ep'$.
Then we can rearrange (\ref{basic}) to obtain for all $0<\r\leq\ep'$ and yet another constant
\begin{align}
\tint_{D_r(y)} e \;\geq\;\tint_{D_{\r r}(\bx)} e 
\;\geq\;  \bc (\r r)^n \bigl( C^{-1}
- a\bc^{\frac 2n} (\r r)^2 - b\bc^{\frac 1n} \r r \bigr).  \label{next}
\end{align}
Now if $a,b\neq 0$ let $\e(a,b)>0$ be the solution of
$$
a \e^2 + b \e = \half C^{-1} .
$$
If $\bc^{\frac 1n} \r r = \e(a,b)$ for some $0<\r\leq\ep'$, then the theorem holds
by the following definition of $\m(a,b)>0$,
\begin{align*}
\tint_{D_r(y)} e \;\geq\; \tfrac 1{2C} \e(a,b)^n \;=:\, \m(a,b).
\end{align*}
Otherwise we must have $\bc^{\frac 1n} \ep' r < \e(a,b)$, so (\ref{next}) with $\r=\ep'$ gives
\begin{align}
\bc \;\leq\; 2C (\ep' r)^{-n} \tint_{D_r(y)} e . \label{last}
\end{align}
In the special case $a=b=0$ we get the same directly from (\ref{next}).
In case $\ep'=\ep$ this proves the theorem since $e(y)\leq 2^n\ep^n\bc$.
Otherwise $\ep'<\ep$ satisfies with another constant
\begin{align*}
2C^{-1} 
&= A_1(\ep' r)^2 + B_1\ep' r + C^{-1} \\
&= \bigl( \sqrt{A_1}\ep'r + C^{-\frac 12} \bigr)^2 + \bigl( B_1 - 2 C^{-\frac 12}\sqrt{A_1} \bigr) \ep'r \\
&= \bigl( \half C^{\frac 12}B_1\ep'r + C^{-\frac 12} \bigr)^2 + \bigl(A_1 - \tfrac C4 B_1^2 \bigr) (\ep'r)^2.
\end{align*}
From this one sees that either $B_1\leq 2C^{-\frac 12}\sqrt{A_1}$ and 
$\ep'r \geq (\sqrt{2}-1)C^{-\frac12} A_1^{-\frac 12}$ from the second line,
or $A_1\leq \frac C4 B_1^2$ and 
$\ep'r \geq 2(\sqrt{2}-1)C^{-1} B_1^{-1}$ from the third line.
Putting this into (\ref{last}) we finally obtain in this last case
$$
e(y) \;\leq\; \bc \;\leq\; C \bigl( A_1^{\frac n2} + B_1^n \bigr) \tint_{D_r(y)} e .
$$

\vspace{-5.5mm}

\hspace*{\fill}$\Box$
\medskip

\subsection*{Acknowledgements}

The mean value inequality with boundary condition originates from my early days as
graduate student, and I am very grateful to Dietmar Salamon for his guidance then
and ever since.
I would also like to thank the participants and organizers of the 2004 Program for Women
in Mathematics at the IAS in \hbox{Princeton} for providing a stimulating atmosphere and
a great test audience for the energy quantization principle.

 \bibliographystyle{alpha}

\end{document}